\documentclass[11pt]{article}

\usepackage{amsmath} 
\usepackage{amsfonts} 
\usepackage{amssymb} 
\usepackage[latin1]{inputenc} 
\newtheorem{theo}{Theorem} 
\newtheorem{Lemma}{Lemma} 
\newcommand\mant{\mathop{\mathcal M}} 
\newcommand\TV{\mathop{\rm TV}} 
\newcommand\R{\mathop{\rm R}}

\begin{document}

\title{The Benford phenomenon for random variables. 
Discussion of Feller's way} 
\date{(version 2, April 19, 2012)} 
\author{Michel {\sc Valadier}\thanks 
{\ Math\'ematiques, Universit\'e Montpellier II, 
place Eug\`ene Bataillon, Case courier 051, 34095 Monptellier Cedex, France. 
email: mivaladier{\string @}wanadoo.fr}} 

\maketitle 

\begin{abstract} 
This is a detailed survey which mainly presents the Pinkham-Feller way. 
I added some new points to the first version \cite{V2} and 
I suppressed ``Examples'' devoted to Gamma, Fréchet and Weibull laws. 
Theorem~\ref{mix} is a bit more general (no assumption of density: this answers a question 
of T.~Hill). 
Section \ref{Poincare} is new and devoted to an argument (Poincaré, Fewster) 
about the effect of high frequencies oscillations. 
Maybe many works, many efforts, have been 
devoted to the study of a sufficient condition of poor value: see Sections 
\ref{catastrophe} and \ref{final}. 
The final Section gives some suggestions. 
\end{abstract} 


\section{Introduction}\label{intro} 
\noindent 
F.~Benford \cite[1938]{B} (and earlier S.~Newcomb \cite[1881]{N}) observed that, in numerical data, 
when the numbers are written in base $10$, very often the first digit, 
which is an integer between $1$ and $9$, takes the value $1$ 
with a frequency much greater than $1/9$ since close to\footnote{\ We 
denote by $\log$ the logarithm in base $10$. 
The logarithm in base $e$ will be denoted by $\ln$.} 
$\log2 = 0.3010..$. 
More generally the \emph{Benford phenomenon} would be that the first digit 
in base $10$, let us denote it by $D$, follows the ``law'': 
\begin{equation}\label{Benf} 
\mathbf P(D = k) = \log\Bigl(\frac{k+1}{k}\Bigr) \quad (k \in \{1,...,9\}) \,. 
\end{equation} 
(Note that $\displaystyle\sum_{k=1}^9 \log\Bigl(\frac{k+1}{k}\Bigr) 
= \log(10) - \log(1) = 1$.) 
\medskip 

\noindent{\sc Comments.} 
The Benford phenomenon is not so intuitive. 
On the contrary: for example integers with $4$ digits go from $1000$ to $9999$. 
``Random'' would give the probability $1/9000$ for each and the 
probability $1/9$ for $\{D = 1\}$. 
For more see Raimi \cite{R1,R2} and Janvresse \cite{J}. 
Not any data can satisfy the Benford phenomenon. 
As said by Scott et Fasli \cite{SF} (just before Section~3) the height of men will give, 
essentially $D = 1$ if expressed in meters, and, if expressed in 
feet\footnote{\ The foot equals $0.3048$ meter.}, 
$D$ will take mainly the values $4$, $5$ and $6$. 
An academic example: if the law of $X$ is uniform on $[1,2]$, $D = 1$ almost surely. 
For negative examples see Section \ref{exneg0}. 
\medskip 

This paper is devoted to the first digit of a random variable and not 
to the other digits or to dynamical systems. 
And I will not discuss papers relying on Fourier Analysis such as 
Pinkham \cite[1961]{P}, Good \cite[1986]{Go}, Boyle \cite[1994]{Bo} 
(note also that Fourier arguments are several times used in \cite{BH2}). 
The litterature is tremendous: Hürlimann \cite{Hu} gives till 
2006, around 350 references; and Berger and Hill \cite{BH3} quote around 600 papers 
(see also \cite{Bee}). 
For French vulgarization papers see \cite{Hi4,D,J}. 
Maybe I missed some important results. 
\medskip 

I will discuss some arguments starting from Feller \cite[1966]{F} and 
quote specially Pinkham \cite[1961]{P}, Engel-Leuenberger \cite[2003]{EL}, 
Dümbgen-Leu\-enberger \cite[2008]{DL}, 
Gauvrit-Delahaye \cite[2008--2009]{GD1,GD2,GD3}, 
Berger \cite[2010]{Br}. 

A mathematical argument going back to Pinkham \cite[1961]{P} is: 
if the density $g$ of $\log X$ is well spread then Benford is approximately satisfied\footnote{\ 
It should be noted that the first part of \cite{P} --- about scale invariance --- is not correct 
(see \cite[Section~4.2 p.40]{BH2}) and that the second part, which is what we quote to, 
relies on Fourier. 
Scale invariance was correctly studied by Hill \cite[1995]{Hi1}, \cite[Th.4.20]{BH2}. 
}. 
In his book Feller \cite[pp.62--63]{F} resumed quickly this result (and quotes Poincaré's roulette) 
with an elementary proof which contains a flaw (see below Section~\ref{Feller}). 
The Feller hypothesis is: $g$ is unimodal and the smallness of the maximum 
ensures the spreadness. 
Exactly the same way can be made correct see \cite{GD1,GD2} (explained in detail in 
Section~\ref{GD} below). 
Moreover Dümbgen and Leuenberger \cite{DL} proved far more better bounds relying firstly on 
total variation of $g$ and further on derivatives of $g$. 
I explained the bounds relying on total variation in Section~\ref{totalv}.

Some ``disaster'' appears: see Section~\ref{catastrophe}. 
Indeed for so usual families of laws on $\mathbb R_+^*$ as the exponential law 
(density $f(x) = \lambda e^{-\lambda \, x}$) or the uniform law (density $(b-a)^{-1}\mathbf1_{[a,b]}$) 
the increasing spreadness of $f$ when $\lambda \rightarrow 0$, resp.\ $b \rightarrow +\infty$, 
is not transmitted to $g$. 
See Section~\ref{??} for examples of the effect of multiplication\footnote{\ As 
shown by \eqref{densite} below, one passes from $f$ to $g$ multiplying essentially by $10^y = x$. 
The maximum, resp.\ the total variation of $g$ equals the maximum (resp.\ total variation) of 
$x \mapsto \ln(10)\, x\,f(x)$.} by $x$, 
\cite{Br} for a carefull discussion of the uniform law 
(resumed in \cite[Prop.4.15 p.\,37]{BH2}), \cite{EL} for the exponential law, 
and \cite{BH1} for a critical review of several arguments (the word ``fallacious'' 
used in this paper, also in \cite[p.\,39 before Th.4.17]{BH2} seems having less pejorative meaning 
in English than ``fallacieux'' has in French). 

Section~\ref{exact} gives naive results. 
Section~\ref{Poincare} analyses some arguments of Fewster's paper \cite{Few} 
going back to Poincaré \cite[1912]{Po}. 
Section~\ref{final} suggests some conclusions. 


\section{Preliminaries}\label{prel} 
\noindent 
Let $X$ be a random variable (briefly r.v.) 
with values in $\mathbb R_+^* = \left]0,+\infty\right[$. 
Let us denote by $D(\omega)$ the first digit in base $10$ of $X(\omega)$. 
It belongs to $\{1,...,9\}$. 
Let $n \in \mathbb Z$ and $k \in \{1,...,9\}$; when $X$ belongs to 
the interval $\left[10^n,10^{n+1}\right[$\,, 
$$ 
D = k \text{ is equivalent to } X \in \left[k \, 10^n,(k+1)10^n\right[ \,. 
$$ 
We abbreviate $\{\omega \,;\, D(\omega) = k\}$ in $\{D = k\}$. 
The following covering is a partition (pairwise disjoint subsets) 
$$ 
\{D = k\} = \bigcup_{n\in\,\mathbb Z} \bigl\{X \in \left[k\,10^n,(k+1)10^n\right[ \bigr\} \,. 
$$ 

The following by Block and Savits \cite{BS} is certainly the most convincing 
early qualitative argument in the direction of Benford: 
if the density $f$ is (strictly) decreasing\footnote{\ Such a density 
has its greatest values near $0$. 
The density $\frac{\mathbf 1_{[10^p,10^q]}(x)}{(q-p)\ln(10)x}$ 
where $p<q$, $p,q \in \mathbb Z$, is null on a neighborhood of $0$ but 
obeys exactly to Benford. 
For more on this, see Remark~1 after Theorem~\ref{th1}.} 
on $\mathbb R_+^* = \left]0,+\infty\right[$ then 
\begin{equation}\label{trivialBS} 
\mathbf P(D = 1) > \mathbf P(D = 2) > \dots > \mathbf P(D = 9) \,. 
\end{equation} 
This comes from the formula 
$$ 
\mathbf P(D=k) = \int_{\bigcup_{n\in\mathbb Z}[k\,10^n,(k+1)10^n[} f(x) \, dx \,. 
$$ 
Note that the gaps in \eqref{trivialBS} can be very small: take $f$ affine on 
$[0,1]$ with a small slope (for example $f(0) = 1 + \varepsilon$, 
$f(1) = 1 - \varepsilon$ and $f(x) = 0$ elsewhere\footnote{\ A strictly decreasing $C^\infty$ 
descent to $0$ when $x$ varies from $1 - \eta$ to $+\infty$ with $>0$ values is possible 
without altering seriously this example.}). 
\medskip 

Plenty of arguments are expressed with the r.v.\ $Y := \log(X)$, 
which takes its values in $\mathbb R$. 
With this r.v.\ the following partition holds 
\begin{equation}\label{I} 
\{D = k\} = \bigcup_{n\in\,\mathbb Z} \bigl\{Y \in \left[\log(k) + n,\log(k+1) + n\right[\bigr\} \,. 
\end{equation} 
Let $\mant(y)$ denotes the \emph{mantissa} of the real number 
$y$ defined by: 
$$ 
\text{if } n\in \mathbb Z \text{ and } y \in \left[n,n+1\right[, \quad \mant(y) := y - n \,. 
$$ 
(The integer $n$ above is the integral part of $y$ usually denoted $\lfloor y \rfloor$ and 
$\mathcal M(y)$ is also called \emph{fractional part} of $y$.) 
Thus $D = k$ is equivalent to (cf.\ \eqref{I}) 
\begin{equation}\label{fond} 
\mant(Y) \in \left[\log k, \log(k+1)\right[ \,. 
\end{equation} 
Assume $Y$ has the density $g$. 
Then $\mant(Y)$ has the density (this is already in \cite[p.\,1224]{P}, \cite [(8.3) p.\,62]{F}, 
\cite[(8.8) p.\,531]{R1}), 
\begin{equation}\label{emp} 
[0,1] \ni y \mapsto \sum_{n \in \mathbb Z} g(n+y) =: \bar g(y) \,, 
\end{equation} 
(the point $1$ should not be in the domain but later for expressing the total 
variation of $\bar g$ it will be useful). 
And let us denote by $\,\overline{\!G}$ the cumulative distribution function of $\mathcal M(Y)$: 
$$ 
\forall z \in [0,1], \quad \,\overline{\!G}(z) = \int_0^z \bar g(u) \, du \,. 
$$ 
We say that $g$ is \emph{unimodal} if $g$ is non-decreasing till some abscissa, and then non-increasing.

The end of this Section is not necessary to understand the remaining of the paper, 
but it corrects the impression caused by the seemingly non-smooth definition of the mantissa. 
Classically the torus $\mathbb T = \mathbb R/\mathbb Z$ is identified to $\left[0,1\right[$. 
With this identification the canonical surjection 
$\varphi : \mathbb R \rightarrow \mathbb T$ coincides with the mantissa $\mathcal M$. 
A geometrical view is: use as for $\mathbb T$ the unit circle $\mathbb U$ via 
the identification 
$$ 
\mathbb T \sim \left[0,1\right[ \ni t \mapsto e^{2\pi i t} \in \mathbb U 
$$ 
and as for $\mathbb R$ the helicoid $\mathbb H$ via the identification 
$$ 
\mathbb R \ni t \mapsto (e^{2\pi i t},t) \in \mathbb H \subset \mathbb U \times \mathbb R\,. 
$$ 
Then $\varphi$ becomes the very smooth map 
$$ 
\mathbb H \ni (z,y) \mapsto z \in \mathbb U \,. 
$$


\section{Position of the problem}\label{Problem} 

\noindent 
We will rarely use the following: 
\begin{equation}\label{BL}
\forall z \in [0,1],\quad \overline{\!G}(z) = z 
\end{equation} 
(equivalently $\bar g(z) = 1$ a.e.) 
in which case one commonly says ``$X$ satisfies the \emph{Benford law}''. 
More correct would be: $X$ modulo $10$ (in the multiplicative group $\mathbb R_+^*$) 
obeys to the Benford law. 

Benford's phenomenon for the first digit is exactly satisfied if 
$$ 
\forall k \in \{2,\dots,9\}, \quad \,\overline{\!G}(\log k) = \log k \,. 
$$ 
As for approximation one can ask for 
$$ 
\forall k \in \{2,\dots,9\}, \quad \,\overline{\!G}\,(\log k) - \log k\ \text{ is small} 
$$ 
or ``$\bar g$ is close to the constant function $\mathbf 1_{[0,1]}$''. 

When we will have found a sufficient condition expressed with $g$ (that is with $Y$) 
the difficulty will be, after expressing it with $X = 10^Y$, to inventory which laws 
satisfy the sufficient condition. 
See Sections~\ref{catastrophe} and \ref{final} where we will question 
about the pertinence of this sufficient condition.


\section{Poincaré's roulette problem (from Feller)}\label{Feller} 

\noindent 
This title comes from Feller \cite[Section~8\,(b) p.\,62]{F}\footnote{\ Note that Feller says 
at the end of his Section~(b) that Poisson's formula could be used. 
And when he turns to Benford in Section~(c) he quotes the name (not any paper) 
of Pinkham.}. 
If a ball is launched from a given zero point on a circle\footnote{\ Poincaré \cite{Po} 
speaks of the roulette pages 11--13, pages 148--150, 
and of digits pages 313--320. 
Two figures in \cite{Few} look a bit like the figure in \cite{Po} page~149. 
See our Section~\ref{Poincare}.} of circumference $1$, 
if the length path is $Y$, the final position of the ball will be $\mathcal M(Y)$. 

When is the law of $\mathcal M(Y)$ close to the uniform distribution? 
Intuitively if one throws the ball with sufficient force and no special effort 
to get an integer number of revolutions or some other precise result, the uniform distribution 
will be approached. 

Feller \cite[p.\,62]{F} says that a valid assumption is ``$g$ is sufficiently spread'' 
(implying a small maximum). 
This is a bit fuzzy. 
As soon after Feller gives a precise and relevant hypothesis: 
the density $g$ of $Y = \log X$ is unimodal and has a small maximum. 
(Below we will quote Pinkham \cite[1961]{P} who worked with a better hypothesis.) 
In my opinion there is a flaw in \cite{F} about which I could not find any precise reference 
in the litterature. 
It is the following: the point $x_k$ defined just after (8.4), 
which is nothing else but 
$\mathcal M(x-a) +a+k$, is not necessarily on the left of $[a+k,a+k+1[$, 
so the assertion, just below (8.5), ``For $k<0$ the integrand is $\leq 0$'' is not correct. 
Nevertheless Raimi \cite[p.\,533]{R1} quotes Feller. 

Under the unimodality hypothesis Gauvrit and Delahaye in 2008 \cite{GD1} gave a correct proof of 
\begin{equation}\label{GDineg} 
\forall z \in [0,1], \quad |\,\overline{\!G}(z) - z| \leq 2 \max g \,. 
\end{equation} 
Their proof is in the line of Feller but they they do not quote him; in \cite{GD2} they spoke 
of ``scatter and regularity'' which are surely not the good words. 
Dümbgen-Leuenberger in 2008 \cite[Th.1 and Cor.2]{DL} still starting from the same 
``spreadness idea'' (they assume that the \emph{total variation} of $g$ is small) 
give far more better bounds: see Section~\ref{totalv} below. 
Already in 1961 Pinkham \cite[Corollary p.\,1229]{P} 
(quoted by Raimi \cite[(8.12) p.\,533]{R1}) using Fourier Analysis arguments obtained 
$$ 
\sup_{0 \leq z \leq 1} |\,\overline{\!G}(z) - z| \leq \TV(g)/6 \,. 
$$ 


\section{The proof of Gauvrit-Delahaye}\label{GD} 

\noindent 
Despite the existence of \cite[1961]{P} and \cite[2008]{DL} 
I reproduce the proof by Gauvrit and Delahaye because it is elementary and pleasant. 
We assume the density $g$ of $Y = \log X$ is unimodal and 
we denote by $M$ its maximum over $\mathbb R$ (maximum to be small). 
\medskip 

\noindent 
\emph{Proof} of \eqref{GDineg}. 
The density $g$ is non-decreasing on 
$\left]-\infty,b\right]$ and non-increasing on $[b,+\infty[$. 
Let $M = g(b)$. 
Without loss of generality we can translate $g$ by an integer $n \in \mathbb Z$, 
so we may suppose\footnote{\ This is the argument in \cite{GD1}. 
Surely $b = 0$ is possible (here we have forgotten $D$ and the factor $10$ relative to $X$) 
and maybe \eqref{GDineg} could be improved of a factor $2$.} $b \in [0,1]$. 
Let $z \in \left]0,1\right]$. 
We will prove the two following inequalities: 
$$ 
\,\overline{\!G}(z) \leq z + 2 \, M \quad \text{and} \quad 
\,\overline{\!G}(z) \geq z - 2 \, M. 
$$ 
The idea (not far from the idea in Feller's book) 
is that on left of $b$ the mean of $g$ on $[n,n+z]$ is less than 
the mean\footnote{\ To prove $\frac{1}{z}\int_0^z g(y) \, dy \leq \int_0^1 g(y) \, dy$ 
when $g$ is non-decreasing on $[0,1]$,
express $\int_0^z g(y) \, dy$ as an integral over $[0,1]$ by a linear change of variable.} 
of $g$ on $[n,n+1]$ and that on right of $b$ the mean of 
$g$ on $[n,n+z]$ is less than that of $g$ on $[n+z-1,n+z]$. 

Precisely: for any $n \leq -1$, since $g$ is non-decreasing on $\left]-\infty,0\right]$, one has 
$$ 
\frac{1}{z} \int_n^{n+z} g(y) \, dy \leq \int_n^{n+1} g(y) \, dy 
$$ 
hence 
\begin{equation}\label{gauche} 
\frac{1}{z} \sum_{n\in\mathbb Z,\,n\leq -1}\int_n^{n+z} g(y) \, dy 
\leq \int_{-\infty}^0 g(y) \, dy \,. 
\end{equation} 
Similarly for any $n \geq 2$ thanks to the non-increasingness of $g$ on $[1+z,+\infty[$, 
$$ 
\frac{1}{z} \int_n^{n+z} g(y) \, dy \leq \int_{n+z-1}^{n+z} g(y) \, dy 
$$ 
hence 
\begin{equation}\label{droite} 
\frac{1}{z} \sum_{n\in\mathbb Z,\,n\geq 2}\int_n^{n+z} g(y) \, dy 
\leq \int_{1+z}^{+\infty} g(y) \, dy \,. 
\end{equation} 
Summing \eqref{gauche} and \eqref{droite} gives 
$$ 
\frac{1}{z} \sum_{n\in\mathbb Z,\,n\neq 0,\,n\neq 1}\int_n^{n+z} g(y) \, dy 
\leq \int_{-\infty}^{+\infty} g(y) \, dy = 1 \,. 
$$ 
On the left hand-side are lacking terms corresponding to $n=0$ and $n=1$. 
Each of them is bounded by 
$$ 
\frac{1}{z} \int_n^{n+z} g(y) \, dy \leq M 
$$ 
hence 
$$ 
\frac{1}{z} \, \,\overline{\!G}(z) \leq 1 + 2 \, M 
$$ 
and 
$$ 
\,\overline{\!G}(z) \leq z + 2 \, M \,. 
$$ 
Now we turn to 
$$ 
\,\overline{\!G}(z) \geq z - 2 \, M \,. 
$$ 
On left of $b$ the mean of $g$ on $[n,n+z]$ is greater than 
the mean of $g$ on $[n+z-1,n+z]$. 
And on right of $b$ the mean of $g$ on $[n,n+z]$ is greater than the mean of $g$ on $[n,n+1]$. 
Thus for $n \leq -1$, 
$$ 
\frac{1}{z} \int_n^{n+z} g(y) \, dy \geq \int_{n+z-1}^{n+z} g(y) \, dy 
$$ 
and for $n \geq 1$, 
$$ 
\frac{1}{z} \int_n^{n+z} g(y) \, dy \geq \int_{n}^{n+1} g(y) \, dy 
$$ 
and summing 
\begin{align*} 
\frac{1}{z} \sum_{n\in\mathbb Z,\,n\neq 0}\int_n^{n+z} g(y) \, dy 
&\geq \int_{-\infty}^{-1+z} g(y) \, dy + \int_1^{+\infty} g(y) \, dy \\ 
&= \int_\mathbb R g(y) \, dy - \int_{-1+z}^1 g(y) \, dy \,. 
\end{align*} 
As the interval $[-1+z,1]$ has length $\leq 2$, the last term 
has absolute value $\leq 2 \, M$. 
$\quad\Box$


\section{Bounds expressed with total variation (Dümbgen-Leuenberger)}\label{totalv} 

We will expose essentially some results by Dümbgen-Leuenberger in 2008 \cite[Th.1 and Cor.2]{DL}. 
Finite total variation encompasses unimodality. 
Precisely if $g$ is unimodal, its total variation is $2\, \max g$. 
With total variation several local minima and maxima are manageable\footnote{\ In \cite{GD1,GD2} 
the authors say that a finite number of bumps is possible. 
The proof could be tedious. 
Example \ref{exneg} below shows that an infinite sequence of bumps may be bad.}. 

Recall that $g$ is the density of $Y$ on $\mathbb R$. 
By the ``stacking'' operation, the density of $\mathcal M(Y)$ on $[0,1]$ is 
$\bar g(z)$ defined in \eqref{emp}. 
A classical notion is \emph{total variation}. 
We assume that $g$ has a finite total variation which we define by\footnote{\ Usually 
this formula is written with strict inequalities. 
It would give the same result (repetition of a value is useless). 
For a fine study of finite total variation functions in one variable, 
but for vector valued functions, see \cite{M}. 
The total variation could be overestimated if one used ``erratic values'' of $g$. 
A non-erratic value at $y$ is a value between the two lateral limits which do exist, 
see for example \cite[Prop.\,4.2 p.\,11]{M}. 
Note that variation is better adapted to cumulative functions than to densities!} 
$$ 
\TV (g) := \sup \Bigl\{\sum_{i=1}^m |g(y_i) - g(y_{i-1})| \,;\, m \geq 1,\ 
-\infty < y_0 \leq \dots \leq y_m < +\infty\Bigr\} \,. 
$$ 
If $g$ is unimodal, $g(y) \rightarrow 0$ when $|y| \rightarrow +\infty$, 
and $\TV(g) = 2 \, \max_\mathbb R g$. 

As for the total variation of $\bar g$ which is a function on the torus $\mathbb T$ 
identified to the half open interval $[0,1[$, 
one should consider 
$$ 
\sup \Bigl\{\sum_{i=1}^m |\bar g(z_i) - \bar g(z_{i-1})| + |\bar g(z_m) - \bar g(z_0)| \,;\, 
m \geq 1,\ 0 \leq z_0 < \dots < z_m < 1 \Bigr\} \,. 
$$ 
But considering $\bar g$ as defined on $[0,1]$ with\footnote{\ Here we can see the importance 
of a sharp definition of $\bar g$: for example if $g(y) = 2 \, y$ on $[0,1]$ and 
$0$ elsewhere, $\bar g(z) = 2 \, z$ on $]0,1[$, but the downfall of $2$ has to be added 
to the progressive increase of $2$ to get the true value $\TV(\bar g) = 4$.} 
$\bar g(1) = \bar g(0)$ one can write 
$$ 
\TV(\bar g) = 
\sup \Bigl\{\sum_{i=1}^m |\bar g(z_i) - \bar g(z_{i-1})| \,;\, 
m \geq 1,\ 0 \leq z_0 \leq \dots \leq z_m \leq 1 \Bigr\} \,. 
$$ 
Now we observe that 
$$ 
\bar g(z) = \lim g_N(z) \quad \text{where} \quad g_N(z) = \sum_{n = -N}^N g(n+z) \,. 
$$ 
Then for any sequence $0 \leq z_0 \leq \dots \leq z_m \leq 1$, 
\begin{align*} 
\sum_{i=1}^m |g_N(z_i) - g_N(z_{i-1})| 
&\leq \sum_{i=1}^m \sum_{n = -N}^N |g(n + z_i) - g(n + z_{i-1})| \\ 
&= \sum_{n = -N}^N \sum_{i=1}^m |g(n + z_i) - g(n + z_{i-1})| \\ 
&\leq \TV(g) 
\end{align*} 
hence the inequality (cf.\ the first assertion of \cite[Theorem~1]{DL} 
and \cite[formula (5) page~107]{DL}), 
$$ 
\TV(\bar g) \leq \TV(g) \,. 
$$ 

Since 
$$ 
\sup \bar g \geq \int_0^1 \bar g(z) \, dz = 1 \geq \inf \bar g 
$$ 
one has 
$$ 
\R(g) := \sup_{z_1 \leq z_2} |\bar g(z_2) - \bar g(z_1)| \geq \sup_z |\bar g(z) - 1| 
$$ 
Note that $|\bar g(z_2) - \bar g(z_1)| = \max([\bar g(z_2) - \bar g(z_1)]^+,[\bar g(z_2) - \bar g(z_1)]^-)$. 
Since $g$ is integrable on $\mathbb R$ it tends to $0$ at infinity, and with the notation 
$$ 
{\TV (g)}^+ := \sup \Bigl\{\sum_{i=1}^m (g(y_i) - g(y_{i-1}))^+ \,;\, m \geq 1,\ 
-\infty < y_0 \leq \dots \leq y_m < +\infty\Bigr\} \,. 
$$ 
and the anologous with negative parts, one has ${\TV}^+(g) = {\TV}^-(g) = \TV(g)/2$. 
Hence (cf.\ \cite[Th.1]{DL}) 
$$ 
\sup_{z\in [0,1]} |\bar g(z) - 1| \leq \R(\bar g) \leq \TV(g)/2 \,. 
$$ 

Now we are going to prove \cite[Cor.2 p.\,102]{DL} 
$$ 
\sup_{0 \leq z_1 < z_2 \leq 1} |(\,\overline{\!G}(z_2) - \,\overline{\!G}(z_1) - (z_2 - z_1)| \leq 
(z_2 - z_1)[1 - (z_2 - z_1)] \, \TV(g)/2 \,. 
$$ 
Let $\delta := z_2 - z_1$. 
Then (we reproduce \cite[proof of Cor.2 p.\,108]{DL}) 
\begin{align*} 
|(\,\overline{\!G}(z_2) - \,\overline{\!G}(z_1) - (z_2 - z_1)| 
&= \Bigl|\int_{z_1}^{z_2} \bar g(z) \, dz - \delta \int_{z_2 - 1}^{z_2} \bar g(z) \, dz\Bigr| \\ 
&= \Bigl|(1-\delta)\int_{z_1}^{z_2} \bar g(z) \, dz - \delta \int_{z_2 - 1}^{z_1} \bar g(z) \, dz\Bigr| \\ 
&= \Bigl|\delta(1-\delta)\int_0^1 \bigl[ \bar g(z_1+\delta t) - \bar g(z_1-(1-\delta) t) \bigr] \, dt \Bigr| \\ 
&\leq \delta(1-\delta)\int_0^1 \bigl| \bar g(z_1+\delta t) - \bar g(z_1-(1-\delta) t) \bigr| \, dt \\ 
&\leq \delta(1-\delta)\R(\bar g)/2 \\ 
&\leq \delta(1-\delta)\TV(g)/2 
\end{align*} 
which implies 
\begin{equation}\label{boundDL} 
\sup_{0 \leq z_1 < z_2 \leq 1} |(\,\overline{\!G}(z_2) - \,\overline{\!G}(z_1) - (z_2 - z_1)| \leq \TV(g)/8 \,. 
\end{equation} 
\medskip 

As already said, in 1961 Pinkham \cite[bottom of page~1228]{P} using Fourier Analysis arguments obtained 
$$ 
\sup_{0 \leq z \leq 1} |\,\overline{\!G}(z) - z| \leq \TV(g)/6 \,. 
$$ 
All these results give better bounds than those of the foregoing Section. 
Indeed, if $g$ is unimodal, \eqref{GDineg} gives 
$$ 
|(\,\overline{\!G}(z_2) - \,\overline{\!G}(z_1) - (z_2 - z_1)| \leq 4 \, \max_\mathbb R g = 2 \, \TV(g) \,. 
$$ 
In their paper \cite{DL} Dümbgen-Leuenberger give other fine bounds when $g$ admits derivatives.


\section{Return to $X$, the disaster}\label{catastrophe} 

\noindent 
Now, what becomes an hypothesis concerning $g$ when expressed in term of $X$ or its density $f$? 
A disaster appears: 
spreadness of $f$ is not equivalent to spreadness of $g$. 
The two reciprocal bijections\footnote{\ The two ordered sets 
$\mathbb R$ and $ \mathbb R_+^*$ are isomorphic.} 
$$ 
\mathbb R \ni y \mapsto 10^y \in \mathbb R_+^* 
\quad\text{and}\quad 
\mathbb R_+^* \ni x \mapsto \log x \in \mathbb R 
$$ 
exchange perfectly the couple $X$ and $Y$ and also the couple of cumulative disstribution 
functions $F_X$ and $F_Y$: one can switch between one and the other 
only by changing $x$ in $10^y$ or $y$ in $\log x$. 
But as for the density one has 
\begin{equation}\label{densite} 
g(y) = \ln(10) \, 10^y \, f(10^y) \quad\text{and}\quad 
f(x) = \frac{g(\log x)}{x \, \ln(10)} \,. 
\end{equation} 
One could switch between the density of $X$ and the density of $Y$ 
only by changing
$x$ in $10^y$ or $y$ in $\log x$ if one had taken for density of $X$ the density of its 
law $\mathbf P_X$ with respect to the following Haar measure\footnote{\ I am indebted 
to J.~Saint-Pierre \cite{SP} for this idea of Haar measure. 
Note that as early as 1970 Hamming \cite{Ha} used the measure 
with density $1/(\ln(10)\,x)$ on the interval $[10^{-1},1]$. 
See also Section \ref{exact}.} on $\mathbb R_+^*$: 
the image (also called push-forward) of Lebesgue on $\mathbb R$ by $y \mapsto 10^y$. 
With respect to the Lebesgue measure this Haar measure has the density 
$x \mapsto [\ln(10)\, x]^{-1}$\,. 

Obviously from \eqref{densite}, unimodality of $g$ is equivalent to unimodality 
of $\tilde f := [x \mapsto x\,f(x)]$ and 
$$ 
\max_{y \in \mathbb R} g(y) = \ln(10) \max_{x \in \mathbb R_+^*} [x \, f(x)]
$$ 
And as for the total variation of $g$, $\TV(g) = \ln(10) \TV(\tilde f)$. 

Here the ``disaster'' occurs: even if $f$ is unimodal, $g$ may be not, see Section~\ref{??}; 
and even if $\max f$ tends to $0$ when a parameter converges to some value, 
the maximum of $x \mapsto x\,f(x)$ may not tend to $0$. 
Despite the fact that log-normal laws (see \ref{LN}) and Pareto laws (see \ref{Pareto}) 
do the work, the uniform law on $[a,b]$ and the exponential law 
(see \ref{unif} and \ref{expo}) exemplify the difficulty. 

As allusively invoked above, classical usual laws described in textbooks are 
families depending on one or several parameters. 
The list is impressive, but the fact that two among the most simple ones fail 
in exemplifying the Benford phenomenon calls for questioning. 
Surely the so many random variables which seem obey to Benford do not follow 
a classical ``usual law'' and the sentence ``if the spread of the r.v.\ 
is very large'' (as in \cite[p.\,63 just after (8.6)]{F}\footnote{\ In Feller 
the r.v.\ is denoted $Y$ but it is the positive variable whose first digit is considered.}) 
is an unwise shortcut. 
For more comments see \cite{Br,BH1} and our Section~\ref{final}. 

The comments in Terence Tao's blog \cite{T1} (see also a book \cite{T2} I did not see) 
are rather informal but their ideas maybe meet those of Hill \cite[Th.3 p.\,361]{Hi1}, 
\cite[Section~6.2 specially Th.6.20]{BH2}. 

As already noticed by many authors, mixing of several data (\cite{Hi1,Hi4,JR} 
and products \cite{Bo} can give good laws. 

I mention that Gauvrit and Delahaye \cite[Th.2]{GD3} say their more general result 
with a strictly increasing function in place of $\log$ applies to more situations.


\section{Examples}\label{ex} 

\subsection{Annoying examples}\label{??} 
One could expect that the hypothesis ``the density $g$ of $Y = \log X$ is unimodal with a 
small maximum'' is usually encountered. 
Expressed with $X$, it means that $x \mapsto x \, f(x)$ is unimodal with a 
small maximum. 
This does not apply to the uniform law and to the exponential law: 
see below Section~\ref{exneg0}. 

Let us give small examples showing the action of multiplication by $x$. 

1) Let 
\begin{equation*} 
f_0(x) = 
\begin{cases} 
x \quad &\text{if $x \in \left]0,1\right]$,} \\ 
x^{-1} \, [1+(x-1)(x-2)/2] \quad &\text{if $x \in [1,2]$,} \\ 
4 \, x^{-3} \quad &\text{if $x \in [2,+\infty[$.} 
\end{cases} 
\end{equation*} 
This is a positive integrable function, so it is, up to a multiplicative coefficient, a density. 
It is decreasing on $[1,2]$ because on this interval 
$$ 
f_0'(x) = \frac{1}{2} - \frac{2}{x^2} \leq 0 
$$ 
so $f_0$ is unimodal. 
But $x \mapsto x\,f_0(x)$ is no longer unimodal. 
It has two maxima, at $x = 1$ and at $x = 2$. 
$\quad\Box$ 

2) Let $f$ defined on $]0,+\infty[$ by $f(x) = 0$ 
on $]0,1/2] \cup \bigcup_{n\geq 1} \{n-1/2\}$, 
$f(n) = 1/n^2$ for all $n \geq 1$ and $f$ affine on all intervals 
$[n-1/2,n]$ and all intervals $[n,n+1/2]$. 
The graph of $f$ consists of a serie of bumps in form of isosceles triangles. 
The total variation is finite with value $2 \sum_{n=1}^\infty 1/n^2$. 
As for $h(x) := x \, f(x)$ this function equals $0$ at 
each $n - 1/2$ and equals $1/n$ at each $n$. 
Since $\sum_{n=1}^\infty 1/n = +\infty$ the total variation of $h$ is infinite. 
$\quad\Box$ 

\subsection{Negative examples}\label{exneg0} 

\subsubsection{A kind of periodicity}\label{exneg} 

If the law of $X$ is carried by the set 
$$ 
\bigcup_{n\in \mathbb Z} [0.9\, 10^n,10^n[ 
$$ 
then $D \overset{\text{a.s.}}{=} 9$ (example inspired by \cite[p.3]{GD2}). 
And this in spite of, as soon as many intervals have $>0$ probabilities, a large ``scattering''. 
This can be realized with a $C^\infty$ density taking strictly positive values 
on each open interval $]0.9\, 10^n,10^n[$. 

\subsubsection{Uniform law}\label{unif} 
The density is $f(x) = \frac{1}{b-a} \, \mathbf 1_{[a,b]}(x)$ (with the parameters 
$a$ and $b$ satisfying $0 \leq a < b$). 
Surely $x \mapsto x\,f(x)$ is unimodal on $\mathbb R_+^*$ and the spread of $f$ 
is large when $b \rightarrow +\infty$. 
The maximum of $x\,f(x)$ is $1$ if $a = 0$. 
And, if $a > 0$ the maximum is $\frac{b}{b-a}$ which decreases when 
$b$ increases to $+\infty$, but the limit is $1$. 
So an inequality as \eqref{GDineg} or \eqref{boundDL} does not apply. 
For a finer study see \cite{Br}, \cite[Prop.4.15 p.\,37]{BH2}. 

\subsubsection{Exponential law}\label{expo} 
The density is $f(x) = \lambda \, e^{-\lambda x}$ ($\lambda \in \left]0,+\infty\right[$ 
is the parameter). 
One could naively expect a good Bendford approximation when $\lambda \rightarrow 0$. 
Derivating $h(x) := x\,f(x)$ one proves easily that 
the function $x\,f(x)$ is unimodal; and its maximum attained at 
$x = 1/\lambda$ has the value $1/e$ (particular case of (10) in \cite{V2} about Gamma law). 
This maximum does not tends to $0$ as $\lambda \rightarrow 0$. 
So an inequality as \eqref{GDineg} does not apply, moreover $4 \, \ln(10) \, e^{-1} = 3.388...$ 
is a very huge value. 
Here \eqref{boundDL} gives 
$$ 
\leq 2\, \ln(10) [\max_{x > 0} x\,f(x)]/8 = \ln(10) \, e^{-1}/4 = 0.211... 
$$ 
Engel and Leuenberger \cite{EL} 
study the exact formula coming from \eqref{fond} 
$$ 
\mathbf P(D=k) = \sum_{n=-\infty}^{+\infty} e^{-\lambda\, k10^n} \, \bigl(1 - e^{-\lambda 10^n}\bigr) \,. 
$$ 
They prove that Bendorf is almost satisfied with a periodical dependance on $\log \lambda$ 
and small gaps. 
But the error does not tend to $0$ as $\lambda \rightarrow 0$. 
Note that from \cite[Fig.1 p.\,363]{EL}, as functions of $\lambda$ the probabilities 
$\mathbf P(D = k )$ oscillates around the ``Benford values'' $\log((k+1)/k)$ 
but they do not take the Benford value simultaneously.

\subsection{Positive examples}\label{expos} 

\subsubsection{Log-normal laws}\label{LN} 

When $X = \exp(Y_0)$ with $Y_0$ of law $\mathcal N(\mu,\sigma^2)$, one has 
$Y = \log X = Y_0/\ln(10)$. Recall that the density of $Y_0$ is 
$$ 
y \mapsto 
\frac{1}{\sigma\,\sqrt{2\pi}} \, \exp\left(-\frac{(y - \mu)^2}{2\sigma^2}\right) \,. 
$$ 
So the required properties of $g$ are clearly verified. 
An inequality as \eqref{GDineg} applies: when $\sigma$ tends to infinity the Benford approximation is good. 

\subsubsection{Pareto laws}\label{Pareto} 

The Pareto law of type~1 (cf.\ \cite[p.7]{GD2}) depends on two 
parameters $\alpha$ and $x_0$ both in $\mathbb R_+^*$ and has the density 
$$ 
f(x) = \frac{\alpha \, x_0^\alpha}{x^{\alpha + 1}} \, \mathbf 1_{[x_0,\infty[}(x) \,. 
$$ 
There $x \mapsto x\,f(x) 
= \displaystyle\frac{\alpha \, x_0^\alpha}{x^\alpha} \, \mathbf 1_{[x_0,\infty[}(x)$ 
is non-decreasing on $\left]-\infty,x_0\right]$ (identically null on $\left]-\infty,x_0\right[$) 
and non-increasing on $[x_0,+\infty[$. 
The maximum reached at $a = x_0$ equals $m = a\,f(a) = \alpha$. 
Thanks to an inequality as \eqref{GDineg}, when $\alpha$ tends to $0$ the probabilities $\mathbf P(D=k)$ 
converge to the values of Benford \eqref{Benf}. 
Note that $X$ has no mean as soon as $\alpha \leq 1$ which indicates a large ``scattering''. 

Pareto laws of type 2 are treated in \cite{GD2}.


\section{Two exact results}\label{exact} 
\noindent 
There exist in the litterature a lot of exact results, some relying on 
``scale invariance'' see \cite{Hi1,Hi2,Hi3}, 
other relying on mixing of laws, see \cite{JR}. 
I will give personal results 
written when I was completely naive with the Benford phenomenon 
and being unaware of \cite{Ha}, \cite{BS} and \cite{BH2}. 

The next Theorem 
shows that the rough hypothesis ``$\mant(\log X)$ follows the uniform law on $[0,1]$'' 
admits sufficient conditions. 
The first part has already been given in 2010 by Block and Savits \cite{BS} 
but there is already in Raimi \cite[1976]{R1} a similar result which would come from Benford: 
see page 532 and the figure taken from Benford; see also \cite[Ex.3.6 p.\,24]{BH2}. 
The second hypothesis comes from the caption of the figure in \cite{GD1,GD2}. 

\begin{theo}\label{th1} 
Let $X$ be a random variable ($X > 0$) and $Y := \log(X)$. 
Suppose that $Y$ has the density $g$. 
Suppose one of the following hypotheses: 

1) $g$ is countably a step function, constant (equality Lebesgue a.e.) 
on each interval $[n,n+1]$ ($n \in \mathbb Z$), i.e.\ 
\begin{equation}\label{step} 
g = \sum_{n \in \mathbb Z} \gamma_n \mathbf 1_{[n,n+1]} \,, 
\end{equation} 

2) $g$ is continuous on $\mathbb R$ and affine on each interval $[n,n+1]$ ($n \in \mathbb Z$). 

\noindent 
Then $X$ follows the Benford law \eqref{BL}. 
\end{theo} 

\noindent 
{\sc Remarks.} 
1) formula \eqref{step} is equivalent to the Block and Savits expression \cite[(3)]{BS} 
(where $\gamma_n$ is $p_n$): 
\begin{equation}\label{BS} 
f(x) = \sum_{n = p}^q \gamma_n \frac{\mathbf 1_{[10^n,10^{n+1}]}(x)}{\ln(10)x} 
\end{equation} 
where $-\infty \leq p \leq q \leq +\infty$, $\gamma_n \geq 0$ and $\sum \gamma_n = 1$. 
Such densities could approximate some real-life densities, but 
a precise study is in the field of Numerical Analysis and Statistics 
(see a discussion in Part~2 of Section~\ref{final}). 
A particular case of \eqref{BS} is the density 
$$ 
f(x) = \frac{\mathbf 1_{[10^p,10^q]}(x)}{(q-p)\ln(10)x} 
$$ 
where $p<q$, $p,q \in \mathbb Z$ which has been given in a footnote of Section~\ref{prel}. 

2) Without the continuity of $g$ the second part does not hold: 
take $g(y) = 2 \, y$ on $[0,1]$ and $0$ elsewhere. 
\bigskip

\noindent 
\emph{Proof of Theorem~\ref{th1}.} 
1) Let $\gamma_n$ be the value of $g$ on $[n,n+1]$. 
The serie \eqref{emp} gives 
$$ 
\sum_{n \in \mathbb Z} g(n+y) = \sum_{n \in \mathbb Z} \gamma_n = 1 
$$ 
which proves that $\mant(\log X)$ follows the uniform law on $[0,1]$. 

2) The integral of $g$ on $[n,n+1]$ is the aera of a trapezoid and it amounts to 
$\frac{1}{2}\,\bigl(g(n)+g(n+1)\bigr)$. 
As the sum is $1$, it holds $\sum_{n = -\infty}^{+\infty} g(n) = 1$. 
Above $[0,1]$ the function 
$$ 
y \mapsto \sum_{n=-N}^N g(n+y) 
$$ 
is affine and equals $\displaystyle\sum_{n=-N}^N g(n)$ at $0$ 
and equals $\displaystyle\sum_{n=-N}^N g(n+1)$ at $1$. 
All this converge to $1$. 
For the incredulous reader if any: the affinity entails 
\begin{align*} 
\sum_{n=-N}^N g(n+y) &= y \sum_{n=-N}^N g(n+1) + (1 - y) \sum_{n=-N}^N g(n) \\ 
&= \sum_{n=-N}^N g(n) + y \, \bigl(g(N+1) - g(-N)\bigr) \\ 
&\rightarrow 1 
\end{align*} 
when $N \rightarrow \infty$. 
$\quad\Box$ 
\medskip 

The next exact result is no longer realistic. 
It as already been obtained by Hamming \cite[Section~{\sc iv} p.\,1615]{Ha} (quoted in 
\cite[p.\,535]{R1}). 
See also \cite[Part (i) of Th.4.13]{BH2} and \cite[Part~1 of Theorem~6.3]{BH2}. 
One could imagine collecting data (richness, level of a river, etc.) 
in several places and several countries where the units are not the same. 
All this would be listed together. 
 
The idea leading to a mathematical result is: multiply a given r.v.\ $X_0$ 
which models our physical quantity (at least in one precise unit) 
by a random coefficient belonging to $[1,10]$, which gives $X$ 
(and as for the law of $X$ a mixing of the laws of the homothetic r.v.\ of $X_0$). 
Changing the unit of several times a factor $10$ or $1/10$ 
would not change the first digit in base $10$. 
We assume that the coefficient obeys the Haar measure\footnote{\ See a foregoing footnote 
in Section \ref{catastrophe}.} 
of the multiplicative group $(\mathbb R_+^*,*)$ 
restricted to $[1,10]$ (more precisely the image of Lebesgue measure by $u \mapsto 10^u$). 
\bigskip

\begin{theo}\label{mix} 
Let $X_0$ be a random variable ($X_0 > 0$) defined on $(\Omega,\mathcal F,\mathbf P_0)$. 
The Lebesgue measure on $[0,1]$ is denoted by $\mathbf\Lambda$. 
Let $X$ be the r.v.\ on $\Omega\times[0,1]$ 
equipped with the probability measure $\mathbf P := \mathbf P_0 \otimes \mathbf\Lambda$ 
defined as 
$$ 
X(\omega,u) = 10^u \, X_0(\omega) 
$$ 
Then $D$ obeys to the Benford law: for $k \in \{1,...,9\}$, 
$\mathbf P(D = k) = \log\bigl(\frac{k+1}{k}\bigr)$. 
\end{theo} 

\noindent 
\emph{Proof.} 
Let $Y_0 := \log(X_0)$. 
For $k \in \{1,...,9\}$ one has 
\begin{align*} 
D(\omega,u) = k &\Longleftrightarrow X(\omega,u) \in 
\bigcup_{n\in\,\mathbb Z} \left[k\,10^n,(k+1)10^n\right[ \\ 
&\Longleftrightarrow u + Y_0(\omega) \in \bigcup_{n\in\,\mathbb Z} \left[n+\log k,n+\log(k+1)\right[ \,. 
\end{align*} 
The above unions are disjoint, hence we have to sum the terms 
\begin{equation*}\label{terme} 
\mathbf P\Bigl(\bigl\{(\omega,u) \,;\, u + Y_0(\omega) \in \left[n+\log k,n+\log(k+1)\right[\bigr\}\Bigr) \,. 
\end{equation*} 
We turn to calculus: 
The transformation of the second line relies on successive integration (firstly 
with respect to $u$ and then to $y$) 
\begin{align*} 
&\mathbf P\Bigl(\bigl\{(\omega,u) \,;\, u + Y_0(\omega) \in \left[n+\log k,n+\log(k+1)\right[\bigr\}\Bigr) \\ 
= &\;(\mathbf P_{Y_0} \otimes \mathbf\Lambda) 
\Bigl(\bigl\{(\omega,u) \,;\, u + y \in \left[n+\log k,n+\log(k+1)\right[\bigr\}\Bigr) \\ 
= &\int_{\mathbb R} \mathbf\Lambda(\left[n+\log k-y,n+\log(k+1)-y\right[ \cap [0,1]) \, d\mathbf P_{Y_0}(y) \\
= &\int_{\mathbb R} \mathbf\Lambda([\log k-y,\log(k+1)-y] \cap [-n,-n+1]) \, d\mathbf P_{Y_0}(y) \,. 
\end{align*} 
But 
\begin{align*} 
\sum_{n\in \mathbb Z} \mathbf\Lambda([\log k-y,\log(k+1)-y] &\cap [-n,-n+1]) \\ 
=\ &\mathbf\Lambda([\log k-y,\log(k+1)-y]) \\ 
=\ &\log(k+1) - \log k \,. 
\end{align*} 
As $\displaystyle\int_{\mathbb R} d\mathbf P_{Y_0}(y) = 1$ this proves 
\begin{multline*} 
\sum_{n=-\infty}^\infty \mathbf P\Bigl(\bigl\{(\omega,u) \,;\, 
Y_0(\omega) \in \left[n-u+\log k,n-u+\log(k+1)\right[\bigr\}\Bigr) \\ 
= \log\Bigl(\frac{k+1}{k}\Bigr) \quad\Box 
\end{multline*} 
\medskip 

\noindent
{\sc Comment.} 
The hypothesis that the unit could be random and obey to a Haar measure is debatable. 
As said by some author, there is a ratio 10 between the decimeter and the meter 
but as for volumes one gets the ratio of 10 between $100$~dm$^3$ and one m$^3$ 
(and not between one dm$^3$ and one m$^3$). 
And usual units are certainly numerous but in a finite number: cf.\ meters and feet 
(argument of \cite{SF} quoted above in Section \ref{intro}). 


\section{The stripey hat of Fewster and Poincaré}\label{Poincare} 
\noindent 
The following result below, convergence \eqref{Few}, is 
with a variant: $n$ in place of $\lambda$, in \cite[Th.4.17 page 39]{BH2}. 
The figure in Fewster \cite[Fig.1 p.\,28]{Few} looks as Poincaré's figure \cite[Fig.15 p.\,149]{Po}. 
In \cite{Few} the grey parts are not half of length but have proportion $z < 1$. 
The idea is to give an intuitive justification of ($g$ being the density of the r.v.\ $Y$) 
$$ 
\int_\mathbb R \bigl(\sum_{n \in \mathbb Z} \mathbf 1_{[n,n+z]}\bigr)(y) \, 
\frac{1}{\lambda}\, g\Bigl(\frac{y}{\lambda}\Bigr) \, dy \longrightarrow z 
\quad\text{as}\quad 
\lambda \rightarrow +\infty \,. 
$$ 
Fewster curve is, as Poincaré's one, ``regular''. 
Moreover the one by Fewster is a bell (or hat) curve. 

The following result is classical for those knowing Young's measures and 
Rademacher functions (cf.\ \cite{V1}). 

\begin{Lemma} 
Let $z \in \left]0,1\right[$ and $\varphi = \sum_{n \in \mathbb Z} \mathbf 1_{[n,n+z]}$. 
Then with the notation 
$$ 
\psi_\lambda(y) := \varphi(\lambda \, y), 
$$ 
the function $\psi_\lambda$ converges, when $\lambda \rightarrow +\infty$, 
$\sigma(L^\infty,L^1)$ to $z\,\mathbf 1_\mathbb R$ that is 
$$ 
\forall h \in L^1(\mathbb R), \quad \int_\mathbb R \psi_\lambda(y) \, h(y) \, dy 
\longrightarrow z \int_\mathbb R h(y) \, dy \,. 
$$ 
\end{Lemma} 

\noindent 
\emph{Proof.} 
The convergence when $h$ is the characteristic function of a compact 
interval, $h = \mathbf 1_{[a,b]}$, is elementary. 
Then it holds for linear combinations of such functions, i.e.\ for functions 
$h$ in a dense subset of $L^1$. 
Since $\psi_\lambda$ belongs to the unit ball of $L^\infty$, an equicontinuous subset 
of the dual space of $L^1$, the result holds for any $h \in L^1$. 
$\Box$ 

Let us denote for $h \in L^1$, $\psi \in L^\infty$, 
$$ 
\langle \psi, h \rangle := \int_\mathbb R \psi(y) \, h(y) \, dy \,. 
$$ 
Now let (note that the action of $\lambda$ is not the same on $\varphi$ and on $g$) 
\begin{equation}\label{dilat} 
g_\lambda(y) = \frac{1}{\lambda}\, g\Bigl(\frac{y}{\lambda}\Bigr) \,. 
\end{equation} 
This is the density of $\lambda\,Y$ (note that $10^{\lambda\,Y} = X^\lambda$; 
and that this dilatation, with $\sigma$ in place of $\lambda$, is already in \cite[page~100]{DL}; 
see also Part~(i) of \cite[Th.6.1]{BH2}). 
By change of variable 
\begin{equation}\label{scal} 
\langle \varphi,g_\lambda\rangle = \langle \psi_\lambda,g\rangle 
\end{equation} 
By \eqref{scal} and Lemma~1
\begin{equation}\label{Few} 
\int_\mathbb R \bigl(\sum_{n \in \mathbb Z} \mathbf 1_{[n,n+z]}\bigr)(y) \, 
\frac{1}{\lambda}\, g\Bigl(\frac{y}{\lambda}\Bigr) \, dy 
= \langle \varphi,g_\lambda\rangle = \langle \psi_\lambda,g\rangle 
\longrightarrow z \int_\mathbb R g(y) \, dy = z \,. 
\end{equation} 
Hence a spreading of $g$ on $\mathbb R$ in the manner of \eqref{dilat} 
(i.e.\ a ``dilatation'') when $\lambda \rightarrow +\infty$ 
implies the expected approximate Benford phenomenon. 
The spreading \eqref{dilat} could be combined with a translation, 
but one is far from the bounds by Dümbgen et Leuenberger \cite{DL}. 
\medskip 

The Poincaré roulette cf.\ \cite[Section~92 pp.148--150]{Po}, \cite{C} 
(and maybe \cite{Fr} which is quoted by \cite{C}, but I did not see it) 
seems an easy result when one knows that the function 
$$ 
r_n(x) = 
\begin{cases} 1 &\text{if $x \in [k/2n,(k+1)/2n[$ for $k$ even,} \\ 
0 &\text{elsewhere} 
\end{cases} 
$$ 
converges $\sigma(L^\infty([0,1]),L^1([0,1]))$ to the constant function $\frac{1}{2}\,\mathbf 1_{[0,1]}$. 
The point is: let a circle be divided in $2n$ arcs of equal lengths, altenatively 
black and red. 
Then if $\bar g$ is a density of probability on the circle, the probability of black 
tends to $1/2$ when $n$ tends to infinity. 
Poincaré assumed that the density is regular. 


\section{Final comments}\label{final} 
\noindent 
1) A sufficient condition may be far from being necessary. 
Surely the reader could find himself many examples. 
I just propose two: 

a) a sufficient condition for a square matrix 
to have a nul determinant is ``the first line is $(0,\dots,0)$''; 

b) a sufficient condition for $x^2$ to be $\leq 4$ ($x$ belonging to $\mathbb R$) 
is ``$-1/2 \leq x \leq 3/2$''; or a worst one: ``$x = -1$''; 
or the tautological ``$x \in \emptyset$''. 

\noindent 
If the proofs are sufficiently involved, if the mathematical objects belong 
to infinite dimensional spaces, it could be hard detecting the ridiculousness of the result. 

Is it true that all work done in the Pinkham-Feller line is of this kind? 

The hope of a good behavior --- with respect to 
Benford --- of a family of laws depending on one or several parameters 
(cf.\ the ``usual laws'' of textbooks) 
when a parameter converges to some limit is surely not the good idea. 
Maybe only some laws (e.g. the log-normal and Pareto's laws (type~I)) perfectly 
realize this hope. 
For some laws (for numerous ones?) the gap is small (see \cite[page 38]{BH2}) 
but convergence to zero does not hold. 
\medskip 

\noindent 
2) Note that as soon as a density expressed by \eqref{BS} is close to the given density $f$, 
the Benford phenomenon will be approximately satisfied. 
How approach a given density $f$ in this way? 

For example which coefficients $\gamma_n$ make \eqref{BS} a valuable approximation 
of\footnote{\ We said more on the exponential law in Section~\ref{expo} 
where some results of Engel-Leuenberger \cite{EL} are outlined.} $f(x) = \lambda\,e^{-\lambda\,x}$? 
Note that on the intervals $[0.01,0.1]$, $[0.1,1]$, $[1,10]$, $e^{-x}$ decreases 
respectively of a factor $1.09$, $2.459$, $8103.08$, 
while $1/x$ decreases of a factor $10$... 
A good adjustement seems difficult. 

In calculus of integrals, integration by the trapezoidal rule seems better than by 
the rectangle method. 
One could appoximate $g$ by trapezoids: the second part of Theorem~\ref{th1} would 
give a density on $\mathbb R$ satisfying Benford. 
But the constraint remains that the intervals are given: they are the $[n,n+1]$. 
\medskip 

\noindent 
3) Limit laws obtained by involved processes could lead to Benford: 
this is suggested by Tao \cite{T1} (see also a book \cite{T2} I did not see) 
and the theorem by Hill \cite[Th.3 p.\,361]{Hi1}, \cite[Section~6.2 specially Th.6.20]{BH2}. 
I reproduce an alinea of \cite[p.\,118]{BH2} about the key hypothesis: 
``Justification of the hypothesis of scale- or base-unbiasedness of significant 
digits in practice is akin to justification of the hypothesis of independence (and 
identical distribution) when applying the Strong Law of Large Numbers or the 
Central Limit Theorem to real-life processes: Neither hypothesis can be formally 
proved, yet in many real-life sampling procedures, they appear to be reasonable 
assumptions.'' 
\bigskip 

\noindent 
{\bf Acknowledgements} 

\noindent 
I wish to thank for different kinds of help H.W.~Block, T.P.~Hill, R.A.~Raimi, T.~de la Rue 
and my colleagues of Nîmes (France) J.-P.~Mandallena and G.~Michaille. 

\providecommand{\bysame}{\leavevmode\hbox to3em{\hrulefill}\thinspace} 
\providecommand{\MR}{\relax\ifhmode\unskip\space\fi MR } 
\providecommand{\MRhref}[2]{%
 \href{http://www.ams.org/mathscinet-getitem?mr=#1}{#2} 
}
\providecommand{\href}[2]{#2}

\end{document}